%% file: Main.tex
\documentclass{article}
\usepackage{amsfonts}
\usepackage{amssymb}
\usepackage{amsthm}
\usepackage{amsmath}
\usepackage{natbib}
\usepackage{hyperref}
\usepackage{url}
\numberwithin{equation}{section}
\usepackage{graphicx}
\usepackage{algorithm}
\usepackage{algorithmic}
\usepackage{enumerate}
\usepackage{boxhandler}
\usepackage{caption}
\usepackage{fancyhdr}

\usepackage{authblk}
\include{defs}

\newcommand{\ud}{\mathrm {d}}
\newcommand{\dInt}{\int}

\renewcommand{\cite}{\citep}

\begin{document}
\title{\Large On the Monotonicity of the Copula Entropy}
\author{Yaniv Tenzer, and Gal Elidan}
\affil{\footnotesize \textit{Department of Statistics, The Hebrew University of Jerusalem}}

\date{}
\pagestyle{myheadings}
\markright{\hspace{0.2in} Y. Tenzer and G. Elidan, On the monotonicity of the copula entropy}

\maketitle

\vspace{-0.4in}
\begin{abstract}
 Understanding the way in which random entities interact is of key interest in numerous scientific fields. This can range from a full characterization of the joint distribution to single scalar summary statistics. In this work we identify a novel relationship between the ubiquitous Shannon's mutual information measure and the central tool for capturing real-valued non-Gaussian distributions, namely the framework of copulas. Specifically, we establish a monotonic relationship between the mutual information and the copula dependence parameter, for a wide range of copula families.  In addition to the theoretical novelty, our result gives rise to highly efficient proxy to the expected likelihood, which in turn allows for scalable model selection (e.g. when learning probabilistic graphical models).

\end{abstract}
\input{Introduction}
\input{Preliminaries}
\input{GeneralTheorem}
\input{ArchimedeanCopulas}
\input{BivariateTwoParameters}

\input{MultivariateArchimedean}
\input{Examples}
\input{Simulations}

\input{Conclusions}

\small{
\bibliographystyle{plainnat}
\bibliography{copula,learn_new}
}
\end{document}

%% file: defs.tex
\long\def\ifundefined#1#2#3{\expandafter\ifx\csname#1\endcsname\relax
#2\else#3\fi}

\ifundefined{DEFNSTYPRESENT}{\def\thmcolon{\hspace{-.75em} {\bf
:} }}{\def\thmcolon{\relax}} 

\newtheorem{THEOREM}{Theorem}[section]
                       {\end{THEOREM}}
\newtheorem{LEMMA}[THEOREM]{Lemma}
\newenvironment{lemma}{\begin{LEMMA} \thmcolon }%
                     {\end{LEMMA}}
\newtheorem{COROLLARY}[THEOREM]{Corollary}
\newenvironment{corollary}{\begin{COROLLARY} \thmcolon }%
                       {\end{COROLLARY}}
\newtheorem{PROPOSITION}[THEOREM]{Proposition}
\newenvironment{proposition}{\begin{PROPOSITION} \thmcolon }%
                           {\end{PROPOSITION}}
\newtheorem{DEFINITION}[THEOREM]{Definition}
\newenvironment{definition}{\begin{DEFINITION} \thmcolon \rm}%
                            {\end{DEFINITION}}
\newtheorem{CLAIM}[THEOREM]{Claim}
                           {\end{CLAIM}}
\newtheorem{OBSERVATION}[THEOREM]{Observation}
\newenvironment{observation}{\begin{OBSERVATION} \thmcolon \rm}%
                           {\end{OBSERVATION}}
\newtheorem{CONJECTURE}[THEOREM]{Conjecture}
                            {\end{CONJECTURE}}
\newtheorem{EXAMPLE}[THEOREM]{Example}
\newenvironment{example}{\begin{EXAMPLE} \thmcolon \rm}%
                            {\end{EXAMPLE}}
\newtheorem{REMARK}[THEOREM]{Remark}
\newenvironment{remark}{\begin{REMARK} \thmcolon \rm}%
                           {\end{REMARK}}
\newtheorem{ASSUMPTION}[THEOREM]{Assumption}
                           {\end{ASSUMPTION}}
\newcommand{\thm}{\begin{THEOREM} \thmcolon \rm}
\newcommand{\thmwname}[1]{\begin{THEOREM} {\bf {#1}} \thmcolon \rm}
\newcommand{\lem}{\begin{lemma}}
\newcommand{\pro}{\begin{proposition}}
\newcommand{\dfn}{\begin{definition}}
\newcommand{\rem}{\begin{remark}}
\newcommand{\xam}{\begin{example}}
\newcommand{\cor}{\begin{corollary}}
\newcommand{\obs}{\begin{observation}}
\newcommand{\prf}{\noindent{\bf Proof:} }
\newcommand{\ethm}{\end{THEOREM}}
\newcommand{\elem}{\end{lemma}}
\newcommand{\epro}{\end{proposition}}
\newcommand{\edfn}{\end{definition}}
\newcommand{\erem}{\bbox\end{remark}}
\newcommand{\exam}{\bbox\end{example}}
\newcommand{\ecor}{\end{corollary}}
\newcommand{\eobs}{\end{observation}}
\newcommand{\eprf}{\bbox}
\newcommand{\beqn}{\begin{equation}}
\newcommand{\eeqn}{\end{equation}}
\newcommand{\beqa}{\begin{eqnarray*}}
\newcommand{\eeqa}{\end{eqnarray*}}

\newcommand{\denselist}{\itemsep 0pt\partopsep 0pt}
\newcommand{\bitem}{\begin{itemize}\denselist}
\newcommand{\eitem}{\end{itemize}}
\newcommand{\benum}{\begin{enumerate}\denselist}
\newcommand{\eenum}{\end{enumerate}}

\newcommand{\commentout}[1]{}
\newcommand{\mynote}[1]{\begin{center}\fbox{\parbox{6in}{\red{#1}}}\end{center}}

\newcommand{\bbox}{\hfill{\vrule height6pt width4pt depth1pt}}

\ifx\boldsymbol\undefined
\newcommand{\boldsymbol}[1]{{\mathbf{#1}}}
\else
\renewcommand{\boldsymbol}[1]{{\mathbf{#1}}}
\fi

\newcommand{\X}{\mathbf{X}}

\ifx\eqref\undefined
\newcommand{\eqref}[1]{Equation~(\ref{#1})}
\else
\renewcommand{\eqref}[1]{Equation~(\ref{#1})}
\fi
\newcommand{\thmref}[1]{Theorem~\ref{#1}}

\newcommand{\figref}[1]{Figure~\ref{#1}}
\newcommand{\secref}[1]{Section~\ref{#1}}

\newcommand{\lemref}[1]{Lemma~\ref{#1}}
\newcommand{\corref}[1]{Corollary~\ref{#1}}
\newcommand{\obsref}[1]{Observation~\ref{#1}}

\def\mymarginpar#1{\mbox{}\marginpar{\raggedleft\footnotesize
  \hspace{0pt}\textsc{#1}}}
\reversemarginpar

\ifx\term\undefined
\newcommand{\term}[3]{\index{#1}\mymarginpar{#2}\nobreak\hspace{0pt}#3}
\else
\renewcommand{\term}[3]{\index{#1}\mymarginpar{#2}\nobreak\hspace{0pt}#3}
\fi

\newcommand{\solution}[1]{}

\newcommand{\answer}[1]{}

\ifx\emcite\undefined
\newcommand{\emcite}[1]{{#1}}
\fi

\newcommand{\bx}{\bf x}

\newcommand{\FX}{F_{\X}(\bx)}
\newcommand{\Fsub}[1]{{#1}}
\newcommand{\Fsubp}[1]{}

\newcommand{\Fxi}[1]{F_{\Fsub{#1}}(x_{#1})}

\newcommand{\fxi}[1]{f_{\Fsub{#1}}(x_{#1})}

\newcommand{\red}[1]{{\textcolor{red}{#1}}}

%% file: Introduction.tex
\section{Introduction}
Understanding the joint behavior of random entities is of great importance in essentially all scientific fields ranging from computational biology and health care to economics and astronomy. Accordingly, the study of joint distributions is fundamental to all the data sciences and goes back at least to the seminal studies of Sir Francis Galton \cite{galton1888co}.

In multivariate modeling, our goal may range from the task of characterizing the full joint distribution which can be difficult, to specific summary statistics such as correlation measures, e.g. Pearson's correlation or mutual information. In this work, we identify a novel and useful relationship between two central frameworks for these tasks, namely the frameworks of copulas and information theory.

In real-valued domains, the most prominent general purpose framework for going beyond the multivariate normal distribution is that of copulas \cite{Joe,Nelsen} pioneered by Sklar \cite{Sklar}.
In a nutshell, copulas allow us to separate the modeling of the (possibly nonparametric) univariate marginals and that of the dependence function. Formally, given a set of univariate marginal cumulative distribution functions $\{F_{X_i}\}$,  a copula function $C_{\bf U}(u_1, \ldots, u_n)$ is a joint distribution over variables that are marginally uniform in the $[0,1]$ range, so that
\[
F_{\bf X}(x_1, \ldots, x_n) = C_{\{F_{X_i}\}}\left(F_{X_1}(x_1), \ldots, F_{X_n}(x_n)\right),
\]
is a valid joint distribution.

This separation between the marginal representation and the copula function that links them allows us, for example, to easily capture multi-modal or heavy-tailed distributions. Indeed, the popularity of copulas as a flexible tool for capturing dependence has grown substantially in recent years \cite{elidan2013copulas}.

Most if not all popular copula families are governed by a dependence parameter that spans the range (or part of it) between independence and full dependence. In fact, in the bivariate case, this is captured by a well known and fundamental relationship between the dependence parameter of the copula and rank-based correlation measures, such as Spearman's rho or Kendall's tau (see, for example, chapter 3 in \cite{Nelsen}). 

\commentout{In this work we establish a relationship between the bivariate copula dependence parameter and %
\[
MI(X;Y) = \int f_{X,Y}(x,y) \log \frac{f_{X,Y}(x,y)}{f_X(x)f_Y(y)} dx dy
\]
or equivalently the \emph{entropy} of the copula density $c_{U \equiv F_X,V \equiv F_Y}(u,v)$
\begin{equation}
\label{Entropy}
    H(C_{U,V}) = \mathbb{E}[-\log c_{U,V}] = -\int c_{U,V}(u,v)\log c_{U,V}(u,v) du dv.
\end{equation}

Importantly, Shannon's mutual information and entropy are fundamental to information theory and are used throughout the exact sciences. Whether in machine learning or physics, it is often
the de facto tool for measuring correlation and/or identifying independence \cite{SOME REFS}. } 

In the field of information theory, Shanon's mutual information measures the reduction in entropy that the knowledge of one variable induces on another. Mutual information is a fundamental tool used to quantify the strength of dependence between random variables \cite{shannon2001mathematical}, and is used throughout
the exact sciences. Formally, Shannon's mutual information is defined as:
\begin{equation}
\label{eq:mi}
MI(X;Y) = \int f_{X,Y}(x,y) \log \frac{f_{X,Y}(x,y)}{f_X(x)f_Y(y)} dx dy.
\end{equation}
Whether in machine learning or physics, mutual information is often the de facto tool for measuring correlation and/or identifying independence \cite{cover2012elements}. A natural question is thus how does this statistic relate to the framework of copulas.

In this work we establish a relationship between the mutual information of two variables and their corresponding copula dependence parameter. Concretely, we prove that the mutual information (equivalently the copula entropy) is monotonic in the bivariate copula dependence parameter for a wide range of copula families, covering the vast majority of copulas used in practice. We also extend our results for the popular class of Archimedean copulas to higher dimensions where other measures of dependence such as Spearman's rho or Kendall tau are not well defined.

The monotonicity result is a theoretical one and an obvious question is whether it has practical merit in the statistical sense. In our simulations we show that it holds in practice using a modest number of samples. An immediate implication is that the mutual information between different pairs of variables can be ranked by evaluating simple statistics that are substantially simpler than the copula entropy, namely Spearman's rho or Kendall's tau. 

Thus, for example, given gene expression data for a large number of genes, we can identify the two genes that have the highest mutual information without ever evaluating the mutual information or the copula density for all pairs of genes, a computationally formidable task. More broadly, our results facilitate highly efficient model selection in scenarios involving many interactions, e.g. when learning probabilistic graphical models. This practical implication is studied in depth in our earlier paper \cite{tenzer2013speedy}, where the theoretical results were substantially more limited.

\commentout{
In this work we aim at establishing the relation between the copula dependence parameter and the \emph{entropy} of the corresponding distribution. Formally, given a continuous density function $f_{\bf X}(x)$, its corresponding entropy is defines as:
\begin{equation}
\label{Entropy}
    H(X) = \mathbb{E}[-\log f(x)] = -\int f(x)\log f(x) dx
\end{equation}
The entropy is a fundamental measure of uncertainty that plays a central role both in statistics and information theory. Concretely, we show that under certain sufficient conditions, the copula entropy function is monotonic in the dependence parameter. This result was previously proven in the case of elliptical families in (complete). However, for other copula families, it was left as a conjecture \cite{elidan2012lightning}.
}

\commentout{Modelling the interaction between random variables, is a fundamental key task, in numerous fields ranging from computational biology to economics to climatology. Naturally, it is the joint distribution that exhibits the dependence structure between the random quantities. Accordingly, in statistics, copulas \cite{Joe,Nelsen} are the central tool for capturing flexible multivariate real-valued distributions by separating the choice of the univariate marginals and the copula function that links them.

Formally, copulas are multivariate distributions whose one-dimensional margins are uniform on the interval $(0,1)$. Alternatively, copulas are functions that join or "couple" multivariate distribution functions to their one-dimensional margins. Concretely, for every multivariate distribution $\FX$ , with a set of marginals $\{F_{x_i}\}$ there exists a copula $C$, such that the following holds:
\[
\FX = C(\Fxi{1},\ldots,\Fxi{n};\theta),
\]
where $\theta$ is the copula \emph{dependence parameter}  that controls the strength of dependence between the variables. \cite{Nelsen, Joe}. By construction, copulas allow one to control both the dependence pattern and the strength of dependence between the variables, separately.   

However, while the joint distribution holds the full picture regarding the dependence structure, it is of a common practice that one is interested only in identifying the strength of dependence between the quantities. 

Accordingly, several \emph{correlation measures} have been designed, some of the most popular are Pearson correlation, Spearman's rho, Kendall's tau and others. These measures are extensively used in every scientific field and serve as summaries of the joint dependence structure. 
}

\commentout{Indeed, the fundamental question of how to quantify dependence between random variables has been the subject of an ongoing research whose scope goes far beyond to that of copulas. Indeed, the mathematical concept of statistical correlation goes back to early days of data analysis by Sir Francis Galton at the end of the 19th century, followed by the seminal work of his prot\`{e}g\`{e} Karl Pearson and the formalization of Pearson's correlation coefficient.}

\commentout{The study of copulas and their application in statistics is a rather modern phenomenon. This growing interest stems out mainly due to the emerging need in flexible statistical modelling. Indeed, copulas are intimately related to the study of distributions with "fixed" or "given" marginals. }

\commentout{From high level view, given a data set, modeling the joint distribution using copulas requires two steps: first is the identification of the dependence pattern between variables and selecting the copula family accordingly. Then, given the selected copula, one need to specify the strength of dependence between the variables, that is, fixing the value of $\theta$. The latter is usually done via a standard maximum-likelihood estimation procedure.} 

\commentout{
A natural question is therefore what is the relation between these measures of association and the copula dependence parameter, $\theta$. 

In the bivariate case for example, it is well known that Spearman's rho between two random variables, say $X$ and $Y$, can be induced directly from their corresponding copula by:
\begin{equation}
\label{Spearmans}
    \rho_{X,Y} = 12\int C(F_X(x),F_Y(y))dxdy - 3.
\end{equation}
A similar relation exists in the case of Kendall's tau correlation \cite{Nelsen}. This kind of results, however, should not be a surprise, as these two measures are based on the ranks, $F_X, F_Y$, whose distribution is fully captured by the copula $C(F_X, F_Y)$. For other measures of dependence though, this question has been left open. In particular for measures that are not ranks based, such as \emph{Pearson correlation}, \emph{Shennon Mutual Information}, etc.
}

\commentout{what is the dependence pattern between the research in this field is conducted along two main axes: the first is a theoretical research that aim at studying copulas fundamental properties and the designing of new copula families with the goal of modelling different dependence structures. The vast majority of works in this category focus on the bivariate case. The second line, however, deals with statistical modeling of high dimensional data using copulas as a building block. Indeed, recently, various high-dimensional constructions that build on a collection of copulas have been suggested. Probably the most notable are the vine \cite{Bedford+Cooke,Kurowicka+Cooke:2002} and the Copula Bayesian Network (CBN) \cite{Kirshner,Elidan:2010}.}

The rest of the paper is organized as follows. In \secref{sec:preliminaries} we briefly review the necessary background on copulas, TP2, super-modular functions and stochastic orders. In \secref{sec:main_result} we present our main theoretical result: broad-coverage sufficient conditions that guarantee the monotonicity of the mutual information (copula entropy) in the copula dependence parameter. We extend the results to the class of bivariate two-parameter families in \secref{sec:2prms}. We then generalize the results for Archimedean copulas of any dimension in \secref{sec:Multv}. In \secref{sec:Simulations} we show the applicability
of our theoretical results in the empirical finite-sample scenario. We finish with concluding remarks in \secref{conclusions}.

%% file: Preliminaries.tex
\section{Preliminaries}
\label{sec:preliminaries}
In this section we briefly review basic definitions and results related to copulas, super modular functions and stochastic orders
that will be needed in the sequel.
\subsection{Copulas}
A copula is a multivariate joint distribution whose univariate marginals are uniformly distributed. Formally: 
\dfn
Let $U_1,\ldots,U_n$ be random variables marginally uniformly
distributed on $[0,1]$. A copula function $C: [0,1]^n \rightarrow
[0,1]$ is a joint distribution
\[
C_\theta(u_1,\ldots,u_n) = P(U_1 \leq u_1,\ldots,U_n \leq u_n),
\]
where $\theta$ are the parameters of the copula function.
\edfn

Now consider an arbitrary set ${\X} = \{ X_1,\ldots X_n \}$ of
real-valued random variables (typically \emph{not} marginally uniformly distributed).
Sklar's seminal theorem \cite{Sklar} states that for \emph{any} joint distribution
$\FX$, there exists a copula function $C$ such that
\[
\FX = C(\Fxi{1},\ldots,\Fxi{n}).
\]
When the univariate marginals are continuous, $C$ is uniquely defined. 

The constructive converse, which is of central interest from a modeling perspective, is also true.  Since $U_i \equiv F_i$ is itself
a random variable that is always uniformly distributed in $[0,1]$,
\emph{any} copula function taking \emph{any} marginal distributions
$\{ U_i \}$ as its arguments, defines a valid joint distribution with
marginals $\{ U_i \}$. Thus, copulas are ``distribution generating''
functions that allow us to separate the choice of the univariate
marginals and that of the dependence structure.

Deriving the joint \emph{density} $f(\bx) = \frac{\partial^nF(x_1,\ldots,x_n)}{\partial
  x_1 \ldots \partial x_n}$ from the copula construction, assuming
$F$ has n-order partial derivatives (true almost everywhere when $F$ is
continuous) is straightforward. Using the
chain rule we can write
\begin{eqnarray}
\label{eq:copula_density}
f(\bx) & = & \frac{\partial^n C(\Fxi{1},\ldots,\Fxi{n})}{\partial \Fxi{1} \ldots \partial \Fxi{n} } 
\prod_i \fxi{i} \\
& \equiv & c(\Fxi{1},\ldots,\Fxi{n}) \prod_i \fxi{i}, \nonumber
\end{eqnarray}
where $c(\Fxi{1},\ldots,\Fxi{n})$, is called the \emph{copula density function}.

Copulas are intimately related to the fundamental concept of \emph{Mutual Information} (MI) of random variables \cite{cover2012elements}. In the bivariate case, the MI of two random variables $X$ and $Y$ is defined as in \eqref{eq:mi}.
Denote $U \equiv F_X$ and $V \equiv F_Y$ so that
$c_{U,V}$ is the copula of the joint distribution 
of $X$ and $Y$ Applying \eqref{eq:copula_density} we get
\begin{eqnarray}
    MI(X;Y) &=&  -\int c_{U,V}(u,v) \log c_{U,V}(u,v) du dv \\ \nonumber &=&  \mathbb{E}[-\log c_{U,V}]  =  H(C_{U,V}).
    \label{MI_Entropy}
\end{eqnarray}
In words, the MI between two random variables equals to the entropy of the corresponding copula. It is easy to see that
this also applies to higher dimensions.

Copulas are also closely tied to other measures of association. The following relationship between the copula function and Spearman's $\rho_{X,Y}$ is well known in the bivariate case:
\begin{equation}
\label{eq:spearman}
    \rho_{X,Y} = 12\int C(F_X(x),F_Y(y))dxdy - 3,
\end{equation}
and a similar relationship is known for Kendall's tau \cite{Nelsen}

\subsection{TP2 and Super-Modular Functions}
Total positive of order two functions (TP2) \cite{olkin2016inequalities} play a central role in statistics and many common copula families have a TP2 density function. Below we define the TP2 concept and provide a simply connection to super-modular functions that will be useful for
our developments.
Formally, denote $u \vee v = min(u,v)$, $u \wedge v = max(u,v).$
\dfn
 A function $\Psi: \mathbb R^2\Rightarrow\mathbb R$ is called TP2 if
\[
  \forall \quad u,v \in\mathbb R^2 \quad  \Psi(u\vee v)\cdot\Psi(u \wedge v) \ge\Psi(u)\cdot \Psi(v).
\]
\edfn
\noindent When the inequality is reversed, the function is reverse rule of order 2 (RR2).
\dfn
A function $\Psi: \mathbb R^2\Rightarrow\mathbb R$ is said to be super-modular if
\[
\forall \quad u,v \in\mathbb R^2,  \Psi(u \vee v)+\Psi(u \wedge v) \ge\Psi(u)+ \Psi(v)
\]
\edfn
\noindent Note that if the inequality is reversed the function is called sub-modular \cite{olkin2016inequalities}.
The following two simple results will be useful in the sequel.


\lem \label{obs:SM}
Let $\Psi(u,v)$ be a positive TP2 (RR2) function. Then $\Phi(u,v)=log(\Psi(u,v))$ is super-modular (sub-modular).
\elem

\prf
From the definition of a TP2 function $\Psi(u \vee v)\cdot \Psi(u
\wedge v)  \ge \Psi(u) \cdot \Psi(v)$, we have 
\beqa
\Phi(u \vee v)+\Phi(u \wedge v) & \equiv & \log(\Psi(u \vee v)) + \log(\Psi(u \wedge v)) \\
  & = & \log(\Psi(u \vee v)\Psi(u \wedge v)) \\
  & \geq & \log( \Psi(u)\Psi(v)) )\\
  & = & \log(\Psi(u)) + \log(\Psi(v)) \\
  & \equiv & \Phi(u) + \Phi(v)
\eeqa
\eprf
%
%
\commentout{
\\
Another result we use is the following: 
\obs \label {obs:TP2} 
Let $f_1(x,y), f_2(x,y)$ be two real non-negative TP2 (RR2) functions. Then $\Psi(x,y)=f_1f_2$ is TP2 (RR2).
\eobs
\prf
Let $(x_1,y_1), (x_2,y_2) \in \mathbb R^2$. Then:
\beqa
\Psi(x_1,y_1)\Psi(x_2,y_2)  & \equiv & f_1(x_1,y_1)f_2(x_1,y_1)f_1(x_2,y_2)f_2(x_2,y_2)\\ 
& = & f_1(x_1,y_1)f_1(x_2,y_2)f_2(x_1,y_1)f_2(x_2,y_2)\\
& \leq & f_1(x_1 \vee x_2,y_1 \vee y_2)f_1(x1 \wedge x_2,y_1\wedge y_2)f_2(x_1 \vee x_2 y_1 \vee y_2))f_2(x_1\wedge x_2),y_1\wedge y_2)\\
& = & f_1(x_1 \vee x_2,y_1 \vee y_2)f_2(x_1 \vee x_2 y_1 \vee y_2))f_1(x1 \wedge x_2,y_1\wedge y_2)f_2(x1\wedge x_2),y_1\wedge y_2)\\
& \equiv & \Psi( x_1 \vee x_2, y_1\vee y_2)\Psi(x_1\wedge x_2,y_1\wedge y_2)
\eeqa
\eprf
}
%
\lem \label {obs:TP2} 
Let $f_1(x,y), f_2(x,y)$ be two real non-negative TP2 (RR2) functions. Then $\Psi(x,y)=f_1f_2$ is TP2 (RR2).
\elem
The proof is immediate and we omit the details. 
\subsection{PQD and Super-Modular Oredrings}
Stochastic orderings introduce the notion of partial orders between random variables. Perhaps the most well known is the  \emph{standard stochastic} order, or using the name more
commonly used in the copula community, the positive quadrant dependent (PQD) order:
\dfn
 Let $\bf{X}, \bf{X'}$ be bivariate random vectors and let $F_{\bf{X}}(u,v),$\\ $F_{\bf{X'}}(u,v)$ be the corresponding distribution functions. $\bf{X'}$ is said to be more PQD then $\bf{X}$ if:
\[
 \forall (u,v) \in \mathbb R^2, F_{\bf{X}}(u,v)\leq F_{\bf{X'}}(u,v)
\]
\edfn
Another stochastic ordering that will be useful in our development is the \emph{Super-modular} (SM) ordering:
\dfn
 Let $\bf{X}, \bf{X'} $ be bivariate random vectors, $\bf{X'}$ is said to be greater than $\bf{X}$ in the super-modular order, denoted by $\bf{X}\leqq_{SM} \bf{X'}, $ if$ \forall\hspace{0.02in} \Psi$ such that $\Psi$  is super modular:
\[
	E_{X}(\Psi(\bf{x}))\leqq E_{X'}(\Psi(\bf{x'}))
\]
\edfn
%
\noindent Note that since most bivariate families are PQD ordered by construction \cite{Joe}, an immediate result of \eqref{eq:spearman} is that Spearman's rho is monotonic in the copula dependence parameter, within a specific copula family.
\commentout{
Another useful result is the following:
\thm \label {Shaked thm}
 ${\bf X}\leq_{pqd} {\bf Y} \iff {\bf X}\leq_{SM} {\bf Y}$ \cite{Shaked}
\ethm
}

%% file: GeneralTheorem.tex
\section{The Monotonicity of Mutual Information in the Dependence Parameter for Bivariate Copulas}
\label{sec:main_result}

We now present our central result and establish a
novel and elegant connection between the copula
function of the distribution of two random variables and the mutual information between these variables (equivalently the copula entropy). In \secref{sec:Multv} we generalize some of these results to higher dimensions.

\thm  \label{thm:Main}
Let $C_\theta(u,v)$ be an absolutely continuous bivariate copula, and let {\bf X}, {\bf X'} be two bivariate random vectors distributed according to the same copula family with two different parameterizations so that  ${\bf X}\sim  C_{\theta_1}(u,v),  {\bf X'} \sim  C_{\theta_2}(u,v)$, where $\theta_2\geq \theta_1\geq 0$. ($\theta_2\leq \theta_1\leq 0$).  Then $-H({\bf X})\leqq-H({\bf X'})$ if one of the following condition holds:
\begin{enumerate}[(a)]
\item $C_\theta(u,v)$ is increasing (decreasing) in $<_{SM}$ and the copula density $c_\theta(u,v)$ is TP2 (RR2) (for all $\theta$).
\label{a}
\item $C_\theta(u,v)$ is an Archimedean copula whose generator $\phi_{\theta}$ is \emph{completely monotone} and satisfies the boundary condition $\phi_{\theta}(0)=1$. In addition the copula is increasing in $<_{SM}$.
\label{b}

\item $C_{\theta,\delta}(u,v)$ is a two-parameters Archimedean copula whose generator $\eta_{\theta,\delta}(s)$ is \emph{completely monotone} and satisfies the boundary condition $\eta_{\theta,\delta}(0)=1$. In addition the copula is increasing in $<_{SM}$, with respect to $\theta$, for a fixed $\delta$.
\label{c}
\item $C_\theta(u,v)$ is an elliptical copula. \label{d}
\end{enumerate}
That is, the negative copula entropy defined in \eqref{MI_Entropy} is monotonic increasing (decreasing) in $\theta$ if any of the above conditions hold for the copula family.
\ethm
We note that in the bivariate case, SM ordering is equivalent to PQD ordering \cite{Shaked} so that the above conditions (a)-(c) can equivalently be stated using a PQD order condition. We also note
that condition (d) was proved in \cite{elidan2012lightning} using an explicit formula of the copula entropy and is stated here for completeness. 

An immediate consequence of \thmref{thm:Main} and the known monotonicity of $\rho_s$ in the dependence parameter $\theta$ for PQD ordered families \cite{Nelsen}, is
\cor
\label{MonotonicityOfSpearman}
If any of the above conditions (a)-(d) hold for a copula family, then the magnitude of Spearman's $\rho_s$ is monotonic in the copula entropy.
\ecor

Importantly, one of the above conditions holds in
most if not all commonly used bivariate copula families so that our result is widely applicable.
In the sections below, we prove the result for each of the sufficient conditions (a)-(c).
In \secref{generalExm}, we survey some example
families.

\subsection{\bf Proof for (a): TP2 Density and SM/PQD Order} To prove the result, we are going to show that the following holds:
\beqa
	\dInt c_{\theta_1}(u,v)log(c_{\theta_1}(u,v))\partial u..\partial u_n & \leqq &  \dInt c_{\theta_2}(u,v)log(c_{\theta_1}(u,v))\partial u\partial v\\  & \leqq & \dInt c_{\theta_2}(u,v)log(c_{\theta_2}(u,v))\partial u\partial v.
\eeqa

Let $\Psi(u,v)=log(c_{\theta_1}(u,v))$. Since $c_{\theta_1}(u,v)$ is TP2, $\Psi(u,v)$  is super modular. ${\bf X}\leqq_{SM}{\bf X'}$. Thus, $E(\Psi({\bf x}))\leqq E(\Psi({\bf x'}))$ for all super modular function $f(u,v)$ and in particular for $f=\Psi(u,v)$, according to ~\lemref{obs:SM}. Thus,  
\beqa
\dInt c_{\theta_1}(u,v)\Psi(u,v) \partial u\partial v \leqq \dInt c_{\theta_2}(u,v)\Psi(u,v) \partial u\partial v.
\eeqa
Now, substituting the explicit form of $\Psi$ we get the first inequality:
\beqa
\dInt c_{\theta_1}(u,v)log(c_{\theta_1}(u,v))\partial u\partial v \leqq \dInt c_{\theta_2}(u,v)log(c_{\theta_1}(u,v)) \partial u\partial v.
\eeqa 

To prove the second inequality we observe that
\beqa
\dInt c_{\theta_2}(u,v)log(c_{\theta_1}(u,v))\partial u\partial v & \leqq & \dInt c_{\theta_2}(u,v)log(c_{\theta_2}(u,v))\partial u\partial v \\ 
& \Leftrightarrow  & \dInt c_{\theta_2}(u,v)log \bigg(\frac{c_{\theta_2}(u,v)}{c_{\theta_1}(u,v)}\bigg)\partial u\partial v \geqq  0 \\
& \Leftrightarrow &  KL(c_{\theta_2} ; c_{\theta_1})  \geqq  0
\eeqa
where $KL(;)$ is the Kullback-Leiber divergence \cite{cover2012elements}, which is always non-negative. Putting
these inequalities together we can conclude:
\begin{eqnarray*}
-H({\bf X}) & \equiv & \dInt c_{\theta_1}(u,v)log(c_{\theta_1}(u,v)) \ud u\ud v  \\ & \leqq & \dInt c_{\theta_2}(u,v)log(c_{\theta_2}(u,v))\ud u\ud v\equiv -H({\bf X'})
\end{eqnarray*}

%% file: ArchimedeanCopulas.tex
\subsection{\bf Proof for (b): Bivariate Archimedean Copulas}
\label{BiArch}
Archimedean copulas are probably the most common non-elliptical copulas that allow for asymmetric
or heavy-tail distributions.

Condition \ref{a} in \thmref{thm:Main} is a fairly general one in the sense that it does not put any restrictions on the copula family type. In particular, it does not require the family to be elliptical or Archimedean. As we shall see below, Condition \ref{b} actually implies condition (a)
in the special case of Archimedean copulas with
completely monotone generator. Among many examples of copula families for which this holds are the Ali-Mikhail (AMH), Clayton, Frank and the Gumbel families (see \secref{generalExm} for details).

\dfn
Let $\psi_{\theta}(x): [0,\infty) \rightarrow [0,\infty]$ be a strictly convex univariate function that is parametrized by $\theta \in \mathbb{R}$. In addition assume that $\psi_{\theta}(\infty)=0$. The Archimedean copula that is generated by $\psi_{\theta}(x)$ is defined as:
\begin{equation} \label{ArchDef}
	C_\theta(u,v) = \psi_{\theta}( \psi_{\theta}^{-1}(u_1)+\psi_{\theta}^{-1}(u_2) ).
\end{equation}
\edfn
\noindent We say that $\psi_{\theta}(x)$ is the copula generator of $C_\theta(u,v)$  \cite{Nelsen}. 
We consider the subclass of copula generators: 
\[
	L_{\infty}=\{ \psi_{\theta}:   (-1)^i\psi_{\theta}(x)^{(i)} \geq 0 \hspace{0.04in} \forall i =0,1,2..,\hspace{0.04in}\psi_{\theta}(0)=1\}.
\]
This is the class of generators whose derivatives alternate signs (this property is widely known as \emph{completely monotonicity}) and in addition these generators satisfy the boundary condition $\psi_{\theta}(0)=1$.

The following lemma shows that an Archimedean copula whose generator $\psi_{\theta}(x)$ is in $L_{\infty}$, can be written as a mixture of two univariate CDFs. We will then use this to show
that \ref{b} implies \ref{a}.
\lem  \label{lem:Mixture}
Let $C_\theta(u,v)$ be an Archimedean copula and let $\psi_{\theta}$ be its generator such that $\psi_{\theta} \in L_{\infty}$. Then there exists a CDF $M(\alpha)$, of a positive random variable $\alpha$,  and unique CDFs $G_1(u), G_2(v)$ such that:
\[
    C_\theta(u,v) = \int_0^{\infty} G_1(u)^{\alpha}\cdot G_2(v)^{\alpha}dM(\alpha).
\]
\elem
For the sack of completeness, we give a simple proof here. Note that a different proof can be found in \citep{Joe}.
\ \\
\prf
From Bernstein's theorem we have that each $\psi_{\theta} \in L_{\infty}$ is a Laplace transform of some distribution function of a positive random variable $\alpha$. That is:
\[
	\psi_{\theta}(s) =\int_0^{\infty} e^{-s\alpha}dM(\alpha), \quad s \geq 0.
\]
Thus
\begin{equation}
\label{eq:c_arch}
C_\theta(u,v)=\psi_{\theta}( \psi_{\theta}^{-1}(u)+\psi_{\theta}^{-1}(v))= \int_0^{\infty} e^{-\alpha( \psi_{\theta}^{-1}(u)+\psi_{\theta}^{-1}(v))}dM(\alpha).
\end{equation}
\hspace{0.04in}In addition, for any arbitrary distribution function $F$, and any positive random variable $M(\alpha)$, there exists a unique distribution function $G$ such that \cite[p.84]{Joe}:
\[
	F(x) = \int_0^{\infty} G(x)^{\alpha} dM(\alpha)=\int_0^{\infty} e^{-\alpha (-\ln G(x))}dM(\alpha)= \psi_{\theta}(-\ln G(x)),
\]
Thus, $G(x)=e^{-\psi_{\theta}^{-1}(F(x))}$. In particular, if $U,V$ are uniform, then $F(u)=u,\hspace{0.03in} F(v)=v$ and therefore there exist $G_1(u), G_2(v)$ such that  $G_1(u)=e^{-\psi_{\theta}^{-1}(u)}$ and $G_2(v)=e^{-\psi_{\theta}^{-1}(v)}$. Substituting this into
\eqref{eq:c_arch}, we get
 \[
	C_{\theta}(u,v)= \int_0^{\infty} G_1(u)^{\alpha}\cdot G_2(v)^{\alpha}dM(\alpha).
\]
\eprf
\ \\
Finally, a mixture representation implies a TP2 density \citep{Joe} . Therefore if in addition PQD/SM ordering holds for these copula families, then condition \ref{a} is implied. 

%% file: BivariateTwoParameters.tex
\subsection{\bf Proof for (d): Bivariate Two-parameters Archimedean Copulas}
\label{sec:2prms}
Two-parameter families are used to capture more than one type or aspect of dependence. For example, one parameter may control the strength of upper-tail dependence while the other indicates concordance or the strength of lower tail dependence. Bivariate Archimedean copulas 
generalize the one-parameter families.
In particular, these families have the
following form:
\beqn  
\label{eq:Bivariate2Prm2}
	C_{\theta,\delta}(u,v)=\psi_{\theta}(-log\hspace{0.02in} K_{\delta}(e^{-\psi_{\theta}(u)},e^{-\psi_{\theta}(v)})),
\eeqn
where $K_{\delta}$ is an Archimedean copula, parametrised by $\delta$, as defined in ~\eqref{ArchDef}, and $\psi_{\theta}$ is a Laplace transform. Let $\phi_{\delta}$ be the generator of $K_{\delta}$. The resulting copula is then also an Archimedean copula with generator $\eta_{\theta,\delta}(s)=\psi_{\theta}(-log \hspace{0.02in} \phi_{\delta}(s))$. For further details see \cite{Joe}. 

For certain families, if we fix the $\delta$ parameter, the copula is then increasing in SM order, with respect to $\theta$. Therefore if the copula generator $\eta_{\theta,\delta}(s)$, is completely monotone and in addition the boundary condition $\eta_{\theta,\delta}(0)=1$ holds, then we are back to the settings of condition ~\ref{c} in \thmref{thm:Main}. As a result we have the following corollary:
\cor 
\label{BivariatePrm}
Let $C_{\theta, \delta}(u,v)$ be a two parameters bivariate copula of the above form specified in ~\ref{eq:Bivariate2Prm2}. Let $\eta_{\theta,\delta}(s)=\psi_{\theta}(-log\hspace{0.02in} \phi_{\delta}(s))$ be its associated generator. Then, if $\eta_{\theta,\delta}(s)$ is completely monotone such that $\eta_{\theta,\delta}(0)=1$ and $C_{\theta,\delta}(u,v)$ is increasing in SM order with respect to $\theta$ for a fixed $\delta$, then the copula entropy is increasing in $\theta$.   
\ecor
The proof is immediate given our previous results and we omit the details. Examples of bivariate two-parameters families for which the conditions given in \corref{BivariatePrm} hold are the BB1, BB2 and BB6 families \cite{Joe}. See \secref{generalExm} for details.

%% file: MultivariateArchimedean.tex
\section{Multivariate Copulas}
\label{sec:Multv}
In this section we generalize some of our results for higher dimensions. Before presenting the main results of this section, let us introduce the following class of positive real univariate functions:
\[
	 L^*_{\infty}=\{ \phi: \mathbb R^{+} \rightarrow \mathbb R^{+} | \phi(0)=0, \hspace{0.02in} \phi(\infty)=\infty,\hspace{0.02in} (-1)^{j-1}\phi^{(j)} \geq 0,\hspace{0.02in} j\geq1 \}.
\]
As we shell see, this class plays a central role in our result regarding multivariate Archimedean copulas. We are now ready to state our main result:  
\thm  \label{thm:MultiMain}
Let $C_\theta(u_1,\ldots, u_n)$ be an absolutely continuous bivariate copula, and let {\bf X}, {\bf X'} be two random n-dimensional vectors such that  ${\bf X}\sim  C_{\theta_1}(u_1,\ldots, u_n),  {\bf X'} \sim  C_{\theta_2}(u_1,\ldots, u_n)$, where $\theta_2\geq \theta_1\geq 0$. ($\theta_2\leq \theta_1\leq 0$).  Then $-H({\bf X})\leqq-H({\bf X'})$ if one of the following conditionד holds:
\begin{enumerate}[(a)]
\item $C_\theta(u_1,\ldots, u_n)$ is increasing (decreasing) in $<_{SM}$ and the copula density $c_\theta(u_1,\ldots, u_n)$ is TP2 (RR2) (for all $\theta$).
\label{multi_a}
\item $C_\theta(u_1,\ldots, u_n)$ is an Archimedean copula whose generator, $\phi_{\theta}$, is \emph{completely monotone} and satisfies the boundary condition $\phi_{\theta}(0)=1$. In addition for all $\theta_1 \leq \theta_2, \quad \phi_{\theta_1}^{-1} \phi_{\theta_2}\in L^*_{\infty}$.
\label{multi_b}
\end{enumerate}
\ethm

It can be easily shown that condition \ref{multi_a} in \thmref{thm:MultiMain} is sufficient in any dimension. We omit the proof since it is essentially identical to the bivariate case. Generalizing condition \ref{multi_b}, however, requires more work which we present in this section.
\subsection{\bf Proof for (b): Multivariate Archimedean Copulas}
Recall that our proof for \ref{b} in the bivariate case relies on the fact that a completely monotone generator that satisfies the boundary condition implies a mixture representation, which in turn implies a TP2 bivariate density. 
We now show that both properties also hold in
the general case of Acrhimedean copulas with completely monotone generators for which the boundary condition holds. To complete the
proof, we will then also need to show that the SM ordering holds under the composition condition $\phi_{\theta_1}^{-1} \phi_{\theta_2}\in L^*_{\infty}$ for $\theta_1 \leq \theta_2$.

The following lemma, proved in \cite[p.85-89]{Joe} substantiates the first property:
\lem  \label{lem:Mixture2}
Let $C_\theta(u_1,\cdot,u_n)$ be Archimedean copula and let $\psi_{\theta}$ be its generator such that $\psi_{\theta}$ is completely monotone and $\psi_{\theta}(0)=1$. Then there exists a CDF $M(\alpha)$ of a positive r.v. and a unique CDF-s $G_1(u_1),\ldots, G_n(u_n)$ such that: $C_\theta(u_1,...,u_n) = \int_0^{\infty} G_1(u)^{\alpha}\ldots G_n(u_n)^{\alpha}dM(\alpha)$
\elem
\noindent The following lemma shows that the second property
also holds:
\lem  \label{lem:MTP2}
Let $C_{\theta}(u_1,\ldots,u_n)$ be a copula that can be represented as in \lemref{lem:Mixture2}. Then $C_{\theta}(u_1,\ldots,u_n)$ has a TP2 density.
\elem

\prf
We assume that $C_{\theta}(u_1,\ldots,u_n)= \int_0^{\infty} G_1(u_1)^{\alpha}\ldots G_n(u_n)^{\alpha}dM(\alpha)$. Taking the derivative with respect to each argument we then have that the copula density is:
\[
    c_{\theta}(u_1,\ldots,u_n)=\alpha^n \prod_i {g_i(u_i) }\int_0^{\infty} G_1(u_1)^{\alpha-1}\cdot..\cdot G_n(u_n)^{\alpha-1}dM(\alpha).
\]
Using $\phi(u_1,\ldots,u_{n-1}, \alpha)=\prod_{i=1}^{n-1}G_i(u_i)^{\alpha-1}$, we
can write
\[
    c_{\theta}(u_1,\ldots,u_n) = \alpha^n \prod_i {g_i(u_i) }\int_0^{\infty} \phi(u_1,\ldots,u_{n-1}, \alpha) \cdot G_n(u_n)^{\alpha-1}dM(\alpha).
\]
Now, for each $i, G_i(u_i)^{\alpha-1}$ is TP2 in $(u_i,\alpha)$. Using \obsref{obs:TP2} we then have that $\phi(u_1,\ldots,u_{n-1}, \alpha)$ is also TP2. From this and \cite{karlin1980classes} we have that $\int_0^{\infty} \phi(u_1,\ldots,u_{n-1}, \alpha) \cdot G_n(u_n)^{\alpha-1}dM(\alpha)$ is TP2 (with respect to $u_1,\ldots,u_{n}$). As $\prod_i {g_i(u_i) }$ is trivially TP2 \cite{Joe}, using \obsref{obs:TP2} again we have that the resulting density is TP2.
\eprf

We have shown that, as in the bivariate case, a completely monotone generator that satisfies the boundary condition implies a TP2 copula density. We now characterize the conditions that also ensure $SM$ ordering. 
\commentout{
Consider the following class of positive real univariate functions:
\[
	 L^*_{\infty}=\{ \phi: \mathbb R^{+} \rightarrow \mathbb R^{+} | \phi(0)=0, \hspace{0.02in} \phi(\infty)=\infty,\hspace{0.02in} (-1)^{j-1}\phi^{(j)} \geq 0,\hspace{0.02in} j\geq1 \}.
\]
}
Concretely $L^*_{\infty}$ provides us with the needed condition via the following theorem \cite{wei2002supermodular}:
\thm 
Let $C_1, C_2$ be two n-dimensional Archimedean copulas and let $\phi_1, \phi_2$ be their associated generators, respectively. If $\phi_1, \phi_2$ are two Laplace transforms (equivalently $\phi_1, \phi_2$ are completely monotone and satisfy the boundary condition $\phi_i(0)=1,\hspace{0.03in} i=1,2$, \cite{Joe}), such that $\phi_1 \phi_2^{-1} \in L^*_{\infty}$, then $C_1 \leq_{SM} C_2$.
\ethm
Putting this and \lemref{lem:MTP2} together we get the desired result.

\commentout{Putting this and \lemref{lem:MTP2} together we get:
\cor
\label{MultMonotonicity}
Let $C_{\theta}(u_1,\ldots,u_n)$ be an Archimedean copula and let $\phi_{\theta}(t)$ be the associated copula generator. If $\phi_{\theta}$ is completely monotone for each $\theta$  and in addition for all $\theta_1 \leq \theta_2, \quad \phi_{\theta_1}^{-1} \phi_{\theta_2}\in L^*_{\infty}$, then its negative entropy is monotonic increasing in the family dependence parameter.
\ecor
}

The conditions of our theorem for multivariate copulas may seem somewhat obscure and therefore not of practical interest. However, we note that they actually apply to some of the most popular multivariate Archimedean copulas, namely the Clayton, Gumbel, Frank and Joe copulas. See \secref{generalExm} for details.

%% file: Examples.tex
\section{Examples}
\label{generalExm}
Below we demonstrate the wide applicability of our theory. We begin with condition \ref{a} of \thmref{thm:Main} that provides the broadest coverage and then also provide examples for the others cases.

\subsection{{\bf Examples Satisfying Condition \ref{a} of \thmref{thm:Main}}}
In all the following examples, the copulas under consideration are known to be increasing in the PQD order and have a TP2/RR2 density function \cite {Joe}. Therefore the entropy for these families is monotonic in $\theta$  by \thmref {thm:Main}:
\begin{itemize}
\item{\bf Bivariate Normal:}
\beqa
  0\leq\theta\leq 1,\quad  C_{\theta}(u,v)=\Phi_{\theta}(\Phi^{-1}(u),\Phi^{-1}(v)), 
\eeqa
where $\Phi$ is the $N(0,1)$ cdf, and $\Phi_{\theta}$ is the BVSN cdf with correlation $\theta$.

\commentout{This copula is increasing in $\leq_c$ order and its density is TP2 \cite {Joe}. Thus its entropy is monotonic in $\theta$  by \thmref {thm:Main}.}
\item{\bf Bivariate Farlie-Gumbel-Morgenstern (FGM):}
\beqa
-1\leq\theta\leq 1, \quad C_{\theta}(u,v)=uv+ \theta uv(1-u)(1-v).
\eeqa
Note that for this family the density is TP2 when $\theta \ge 0$ and RR2 when $\theta \leq 0$ \cite {Joe}. Thus its negative entropy is monotonic increasing in $\theta$, when $\theta \in (0,1]$ and monotonic decreasing when $\theta \in [-1,0)$ \cite{Joe}. Therefore overall the FGM copula negative entropy is monotonic increasing in $|\theta|$.
\commentout{This copula is increasing in $\leq_c$ order and its density is TP2 when $\theta \ge 0$ and RR2 when $\theta \leq 0$ \cite {Joe},. Thus its negative entropy is monotonic increasing in $\theta$, when $\theta \in (0,1]$ and monotonic decreasing when $\theta \in [-1,0)$. In other words, the FGM copula negative entropy is monotonic increasing in $|\theta|$.}
\item{\bf Bivariate Frank:}
\beqa
0\leq\theta\leq\infty, \quad C_{\theta}(u,v)=-\theta^{-1}log\bigg(\big[\tau-\big(1-e^{-\theta u}\big)(1-e^{-\theta v}\big)\big]/\tau\bigg),
\eeqa
where $\tau=1-e^{-\theta}.$  
\commentout{This copula is increasing in $\leq_c$ order and its density is TP2 \cite {Joe}. Thus its negative entropy is monotonic increasing in $\theta$. }
\item{\bf Bivariate Gumbel:}
\beqa
1\leq\theta\leq\infty, \quad C_{\theta}(u,v)=e^{-[(\hat u)^\theta+(\hat v)^\theta)]^{1/\theta}},
\eeqa
where $\hat u=-log(u), \hat v=-log(v)$. 
\commentout{This copula is increasing in $\leq_c$ order and its density is TP2 \cite {Joe}. Thus its negative entropy is monotonic increasing in $\theta$.}
\item{\bf Bivariate Clayton:}
\beqa
0\leq\theta\leq\infty, \quad C_{\theta}(u,v)=(u^{-\theta}+v^{-\theta}-1)^{(-1/ \theta)}. 
\eeqa
\commentout{This copula is increasing in $\leq_c$ order and its density is TP2 \cite {Joe}. Thus its negative entropy is monotonic increasing in $\theta$.} 
\commentout{\item{\bf Bivariate Joe:}
\beqa
1\leq\theta\leq\infty, \quad C_{\theta}(u,v)=1-[\Bar{u}^\theta +\Bar{v}^\theta-\Bar{u}^\theta\Bar{v}^\theta]^{1/\theta},
\eeqa
where $\Bar{u}=1-u, \Bar{v}=1-v$. 
\commentout{This copula is increasing in $\leq_c$ order and its density is TP2 \cite {Joe}. Thus its negative entropy is monotonic increasing in $\theta$.}}
\end{itemize}
\subsection{{\bf Archimedean Copula Examples Statisfying condition \ref{b} in \thmref{thm:Main}}}
In the following we demonstrate condition \ref{b} for several Archimedean families that are not covered by
the previous section such as the Clayton/Frank/Gumbel copulas. These examples are slightly lesser known but still useful Archimedean copula families. In all the following examples, the copulas under consideration are known to be increasing in the PQD order \cite{Joe}. In addition it can be easily verified that the boundary condition $\psi_{\theta}(0)=1$ holds.  

In order to establish the completely monotonicity (and hence also the total positivity of the corresponding density function) of the copula generator, we use the following sufficient conditions \cite{Nelsen, widder1942completely}:
\lem 
Let $f(x), g(x)$ be two real univariate functions and let $h_1(x)=f(x) \circ g(x),\hspace{0.03in} h_2(x)=f(x)g(x)$. Then: 
\begin{enumerate}[(i)]
\item If g is completely monotonic and f is absolutely monotonic, i.e., $\frac{\partial f(x)} {\partial x^k} \geq 0$ for k = 0,1,2,... then $h_1(x)$ is completely monotone.

\item If f is completely monotonic and g is a positive function with a
completely monotone derivative, then $h_1(x)$ is completely monotone.
\item If f and g are completely monotone, then so is $h_2(x)$.
\end{enumerate}
\elem
\begin{itemize}

\item {\bf Bivariate Ali-Mikhail (AMH):}
\beqa
\theta \in [-1, 1], \quad	C_{\theta}(u,v)= \frac{uv}{1-\theta(1-u)(1-v)}, 
\eeqa
The generator of this copula is given by:
\[
	\psi_{\theta}(t)=\frac {1-\theta}{exp(t)-\theta}.
\]
This generator is completely monotone for $\theta \in (0,1]$ \citep{jaworski2010copula}. When $\theta \in [-1,0)$, it can easily been shown that the corresponding densities are RR2, using  \lemref{obs:TP2}. Therefore overall the AMH copula negative entropy is monotonic increasing in $|\theta|$.
\item  {\bf Bivariate Joe:}
\beqa
\theta \in [1,\infty), \quad	C_{\theta}(u,v)= 1-\left[(1-u)^{\theta}+(1-v)^{\theta}-(1-u)^{\theta}(1-v)^{\theta}  \right]^{1/\theta}, 
\eeqa
The generator of this copula is given by:
\[
	\psi_{\theta}(t)=1-(1-\exp(-t))^{1/\theta}.
\]
Taking $f(t)=1-t^{1/{\theta}}$ and $g(t)=1-\exp(-t)$, we see that $f$ is completely monotone and $g$ is a positive function whose first derivative completely monotone. Thus their composite is also completely monotone. 

\item  {{\bf Family 4.14} \cite {Nelsen}:}
\beqa
	\theta \in [1,\infty), \quad C_{\theta}(u,v)= \left(1+ \left( (u^{-1/\theta}-1)^{\theta}+(v^{-1/\theta}-1)^{\theta}\right)^{\theta} \right)^{-\theta},  
\eeqa
The generator of this copula is given by:
\[
	\psi_{\theta}(t)=(t^{1/\theta}+1)^{-\theta}.
\]
By taking $f=(1+t)^{-\theta},\hspace{0.03in} g=t^{1/\theta}$ and repeating the same arguments as in previous examples, we get that this generator is completely monotone.
\item  {{\bf Family 4.19} \cite {Nelsen}:}
\beqa
\theta \in (0,\infty), \quad	C_{\theta}(u,v)= \frac {\theta}{\ln (e^{\theta/u} + e^{\theta/v} - e^{\theta})}, 
\eeqa
The generator of this copula is given by:
\[
	\psi{\theta}(t)=\frac {\theta}{\ln \left(t+\exp(\theta)\right)}.
\]
By taking $f=\theta/t,\hspace{0.03in} g= \ln(t+\exp(\theta))$ and repeating the same arguments as in previous examples, we get that this generator is completely monotone.
\commentout{
\item {{\bf Family 4.20} \cite {Nelsen}:}
\beqa
\theta \in (0,\infty), \quad C_{\theta}(u,v)= \left(\ln\left( exp(u^{-\theta})+ exp(v^{-\theta})- e \right)\right)^{-1/\theta},  
\eeqa
The generator of this copula is given by:
\[
	\psi_{\theta}(t)=(\ln(t+e))^{-1/\theta}.
\]
Taking $f=x^{1/\theta}$ and $g=log(t+e)$, we see that this generator is completely monotone since $f$ is completely monotone and $g$ has a completely monotone first derivative.}
\end{itemize}
\subsection{{\bf Examples of Bivariate Two-parameter Families}}

We now provide some examples of two-parameter copula families that satisfy condition \ref{c} of \thmref{thm:Main}. That is, they are all positively PQD/SM ordered with respect to $\theta$ and have a completely monotone generator that satisfies the boundary condition \cite{Joe}. Therefore by ~\corref{BivariatePrm} for fixed $\delta$, their entropy is monotonic in $\theta$.

Recall that bivariate two-paraneter families have the
following form:
\beqn  
\label{eq:Bivariate2Prm}
	C_{\theta,\delta}(u,v)=\psi(-log\hspace{0.02in} K_{\delta}(e^{-\psi_{\theta}(u)},e^{-\psi_{\theta}(v)})),
\eeqn
where $K_{\delta}$ is an Archimedean copula, parametrised by $\delta$, as defined in ~\eqref{ArchDef}, and $\psi_{\theta}$ is a Laplace transform. Let $\phi_{\delta}$ be the generator of $K_{\delta}$. The resulting copula is then also an Archimedean copula with generator $\eta_{\theta,\delta}(s)=\psi_{\theta}(-log \hspace{0.02in} \phi_{\delta}(s))$. For further details see \cite{Joe}.

\begin{itemize}
\item {\bf Family BB1 \cite {Joe}:} Taking K to be Gumbel copula and $\psi_{\theta}(s)=(1+s)^{-1/\theta}, \quad \theta \geq 0$. The resulting copula is:
\[
    C_{\theta,\delta}(u,v)=\left(1+\left((u^{-\theta}-1)^{\delta}+(v^{-\theta}-1)^{\delta}\right)^{1/\delta} \right)^{-1/ \theta},
\]

where $\eta_{\theta,\delta}(s)=(1+s^{1/ \delta})^{-1/\theta}$. 
\commentout{This generator is completely monotone \cite{Joe}. In addition this family is positively PQD ordered with respect to $\delta$, and thus is also SM ordered with respect to $\delta$ \cite{Joe}. Thus, by ~\corref{BivariatePrm} for fixed $\theta$, the copula entropy is monotonic in $\delta$.}

\item {\bf\cite {Joe}, bivariate BB2:} Taking K to be Clayton copula and $\psi_{\theta}(s)=(1+s)^{-1/ \theta}, \quad \theta \geq 0$ the resulting copula is:
\[
    C_{\theta,\delta}(u,v)=\left( 1+ \delta^{-1}log \left ( e^{\delta(u^{-\theta}-1)}+e^{\delta(v^{-\theta}-1)}\right) -1\right)^{-1/\theta},
 \] 

where $\eta_{\theta,\delta}(s)=[1+\delta^{-1}log(1+s)]^{-1/\theta}$. 
\commentout{Repeating the same arguments as in previous example, it can easily been shown that this generator is completely monotone. In addition this family is positively PQD ordered with respect to $\delta$ and thus is also SM ordered with respect to $\delta$ \cite{Joe}. Hence by \corref{BivariatePrm} for fixed $\theta$, the copula entropy is monotonic in $\delta$.}

\item {\bf\cite {Joe}, bivariate BB6:} Taking K to be Gumbel copula and $\psi_{\theta}(s)=1-(1-e^{-s})^{1/ \theta}, \quad \theta \geq 1$ the resulting copula is:
\[
    C_{\theta,\delta}(u,v)=1 - \left ( 1-exp\left(-\left((-\log(1-\hat u ^{\theta}))^{\delta}+(-\log(1-\hat v^{\theta}))^{\delta}\right)^{\frac{1}{\delta}}\right) \right ) ^{\frac {1}{\theta}},
 \]

where $\hat u=1-u, \hat v= 1-v,\quad \eta_{\theta,\delta}(s)=1-[1-\exp\left(-s^{1/\delta}\right)]^{1/\theta}$.
\commentout{Repeating the same arguments this generator it cab easily been shown that this generator is completely monotone \cite{Joe}. In addition this family is positively PQD ordered with respect to $\delta$ and thus is also SM ordered with respect to $\delta$ \cite{Joe}. Hence by ~\corref{BivariatePrm} for fixed $\theta$, the copula entropy is monotonic in $\delta$.}
\end{itemize}
\subsection{{\bf Examples of Multivariate Copulas}}

We finish with some multivariate examples that satisfy the conditions of \thmref{thm:MultiMain}. That is, in each example the copula generator, $\phi_{\theta}$, is completely monotone and in addition $\phi_{\theta_1}^{-1} \phi_{\theta_2} \in  L^*_{\infty}$ for $\theta_1 \leq \theta_2$ \cite{Joe}.
\begin {itemize}

\item {{\bf Multivariate Clayton:}}
\beqa
\theta_1 \leq \theta_2 \in (0,\infty), \quad \phi_{\theta_1}^{-1} \phi_{\theta_2}=\frac{1}{\theta_1}(\theta_2+1)^{\theta_1/ \theta_2} -\frac{1}{\theta_1}.
\eeqa
\commentout{It can be easily shown that $\phi_{\theta_1}^{-1} \phi_{\theta_2} \in  L^*_{\infty}$.
\mynote{Did we actually show it or just cite? If not, better not to say "we have already shown". Also applies to the following example}
In addition, we have already shown that $\phi_{\theta}$ is completely monotone for all $\theta$. Thus the negative entropy of the multivariate Clayton is monotonic increasing in $\theta$.}

\item{{\bf Multivariate Gumbel:}}
\beqa
\theta_1 \leq \theta_2 \in (1,\infty), \quad \phi_{\theta_1}^{-1} \phi_{\theta_2}=t^{\theta_1/ \theta_2}.
\eeqa
\commentout{It can be easily shown that $\phi_{\theta_1}^{-1} \phi_{\theta_2} \in  L^*_{\infty}$. In addition, we have already shown that $\phi_{\theta}$ is completely monotone for all $\theta$. Thus, the negative entropy of the multivariate Gumbel is monotonic increasing in $\theta$.}

\item{{\bf Multivariate Frank:}}
\beqa
\theta_1 \leq \theta_2 \in (0,\infty), \quad \phi_{\theta_1}^{-1} \phi_{\theta_2}=- \ln \left(\frac{(e^{-t}(e^{-\theta_2}-1)+1)^{\theta_1/ \theta_2}-1}{e^{-\theta_1}-1}\right).
\eeqa
\commentout{It can be easily shown that $\phi_{\theta_1}^{-1} \phi_{\theta_2} \in  L^*_{\infty}$. In addition, we have
already shown that $\phi_{\theta}$ is completely monotone for all $\theta$. Thus, the negative entropy of the multivariate Frank is monotonic increasing in $\theta$.}

\commentout{\item{{\bf Multivariate Ali-Mikhail:}} 
\beqa
\theta_1 \leq \theta_2 \in (0,1), \quad \phi_{\theta_1}^{-1} \phi_{\theta_2}=\ln \frac{e^{t}(1-\theta_1)+ \theta_1-\theta_2}{1-\theta_2}.
\eeqa
}
\commentout{It can be easily shown that $\phi_{\theta_1}^{-1} \phi_{\theta_2} \in  L^*_{\infty}$. In addition, we have already that $\phi_{\theta}$ is completely monotone for all $\theta$. Thus, the negative entropy of the multivariate Ali-Mikhail is monotonic increasing in $\theta$.}

\item{\bf Multivariate Joe :}
\beqa
\theta_1 \leq \theta_2 \in (0,1), \quad \phi_{\theta_1}^{-1} \phi_{\theta_2}=-\ln \left(1-(1-\exp(-t))^{\theta_1/\theta2} \right).
\eeqa
\end{itemize}

%% file: Simulations.tex
\section{Simulations}
\label{sec:Simulations}
 \begin{figure*}
 \captionsetup{width=1.3\linewidth}
 \hspace{-1in}
 \begin{minipage}{0.700\columnwidth}
\begin{center}
\begin{tabular}{c}
\begin{tabular}{cc}
\begin{minipage}{0.700\columnwidth}
\hspace{-0.35in}
   \vspace{0.3in} {\large $d=2$}
\end{minipage}&
\hspace{-2.4in}
\begin{tabular}{cccc}
  Clayton & Frank & Gumbel & Joe \\
  \includegraphics[trim=0.2in 0in 0.65in 0.65in,clip,height=0.45\columnwidth]{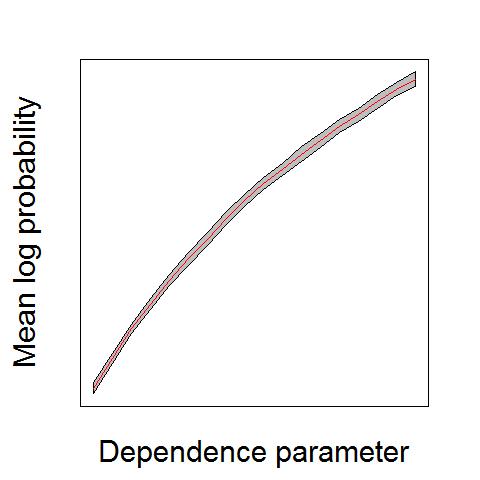} & \hspace{-0.17in}\includegraphics[trim=0.2in 0in 0.65in 0.65in,clip,height=0.45\columnwidth]{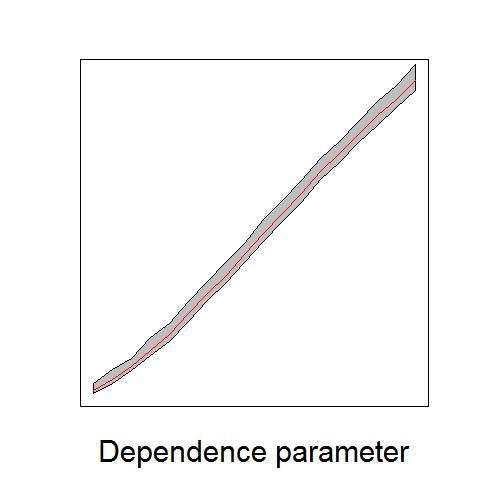} & \hspace{-0.2in} \includegraphics[trim=0.2in 0in 0.65in 0.65in,clip,height=0.45\columnwidth]{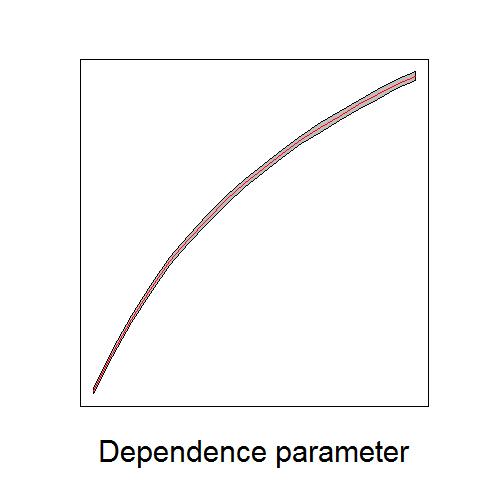}& \hspace{-0.17in}\includegraphics[trim=0.2in 0in 0.65in 0.65in,clip,height=0.45\columnwidth]{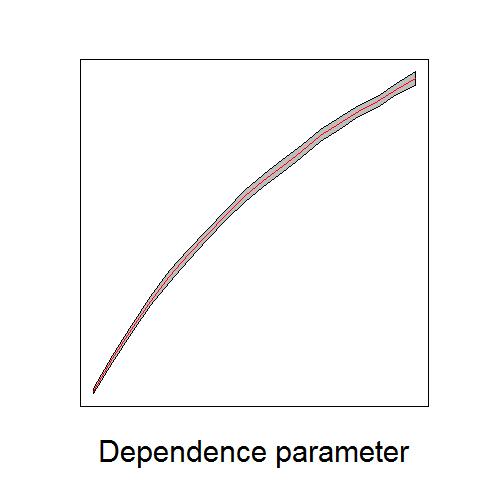}\vspace{0.2in}
\end{tabular}
\end{tabular}\\
\begin{tabular}{cc}
\begin{minipage}{0.700\columnwidth}
\hspace{-0.35in}
 \vspace{0.3in} 
    {\large $d=2$}
\end{minipage}&
\hspace{-2.4in}
\begin{tabular}{cccc}
 Gaussian & T, $\nu=3$ & T, $\nu=7$ & AMH \\
\includegraphics[trim=0.2in 0in 0.65in 0.65in,clip,height=0.45\columnwidth]{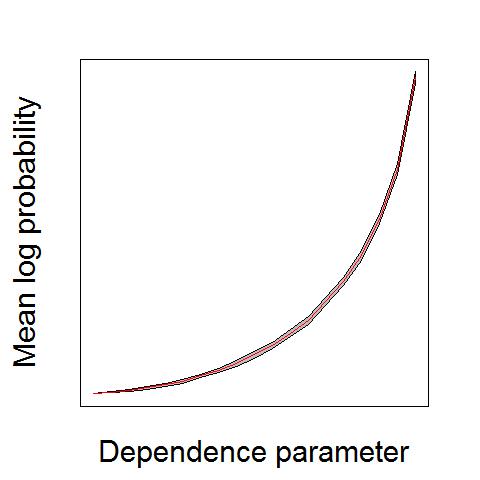} & \hspace{-0.17in}\includegraphics[trim=0.2in 0in 0.65in 0.65in,clip,height=0.45\columnwidth]{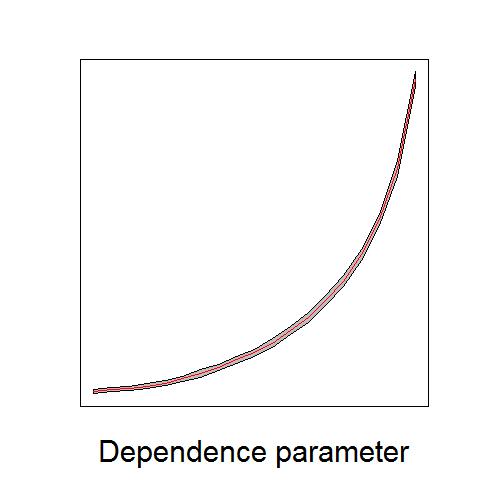} & \hspace{-0.2in} \includegraphics[trim=0.2in 0in 0.65in 0.65in,clip,height=0.45\columnwidth]{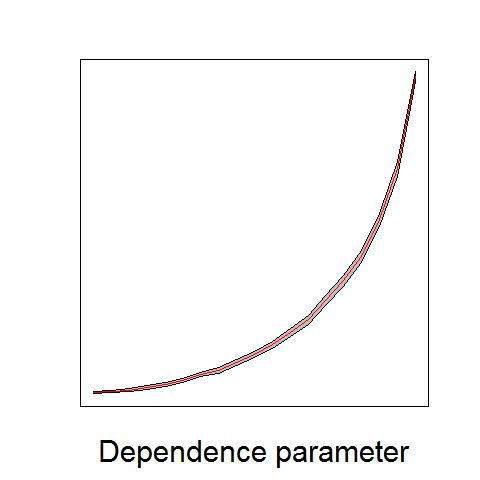}&\hspace{-0.17in} \includegraphics[trim=0.2in 0in 0.65in 0.65in,clip,height=0.45\columnwidth]{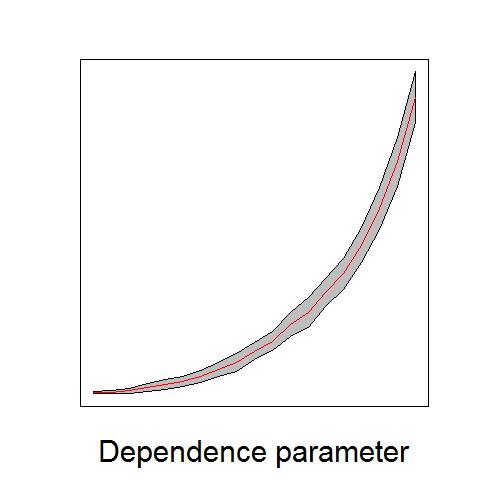}\vspace{0.2in}
\end{tabular}
\end{tabular}\\

 \begin{tabular}{cc}
\begin{minipage}{0.700\columnwidth}
\hspace{-0.35in}
 \vspace{0.3in} 
    {\large $d=2$}
\end{minipage}&
\hspace{-2.4in}
\begin{tabular}{cccc}
 BB1, $\delta=1.5$ & BB1, $\delta=3$ & BB6, $\delta=1.5$ & BB6, $\delta=3$ \\
 \includegraphics[trim=0.2in 0in 0.65in 0.65in,clip,height=0.45\columnwidth]{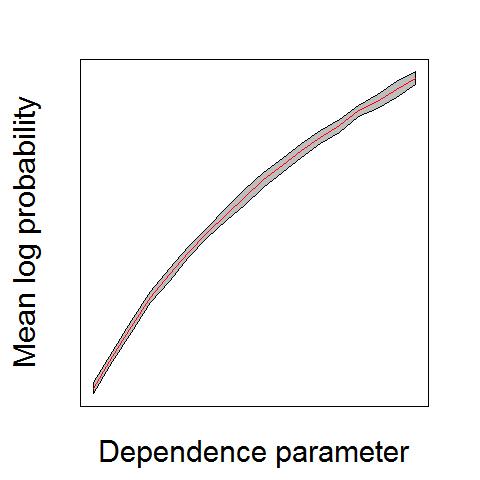} & \hspace{-0.17in}\includegraphics[trim=0.2in 0in 0.65in 0.65in,clip,height=0.45\columnwidth]{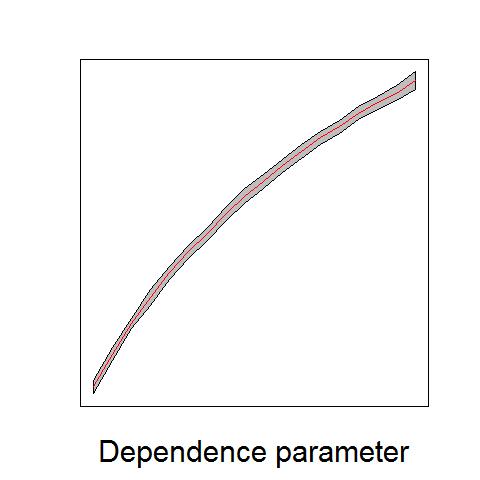} & \hspace{-0.2in} \includegraphics[trim=0.2in 0in 0.65in 0.65in,clip,height=0.45\columnwidth]{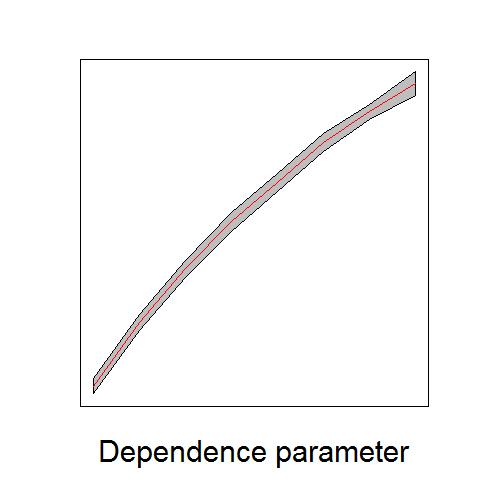}& \hspace{-0.17in}\includegraphics[trim=0.2in 0in 0.65in 0.65in,clip,height=0.45\columnwidth]{Figures/BB1_theta_4.jpg}\vspace{0.2in}
 \end{tabular}
\end{tabular}\\
\begin{tabular}{cc}
\begin{minipage}{0.700\columnwidth}
\hspace{-0.35in}
 \vspace{0.3in} 
    {\large $d=5$}
\end{minipage}&
\hspace{-2.4in}
\begin{tabular}{cccc}
Clayton & Frank & Gumbel & Joe \\
 \includegraphics[trim=0.2in 0in 0.65in 0.65in,clip,height=0.45\columnwidth]{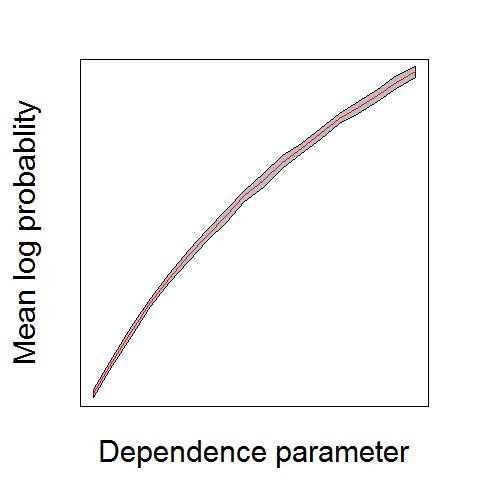} & \hspace{-0.17in}\includegraphics[trim=0.2in 0in 0.65in 0.65in,clip,height=0.45\columnwidth]{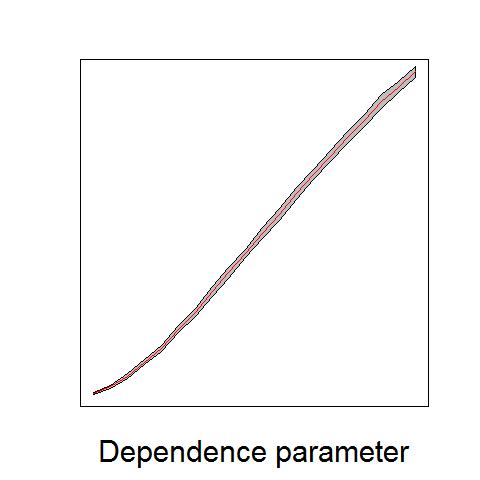} & \hspace{-0.2in} \includegraphics[trim=0.2in 0in 0.65in 0.65in,clip,height=0.45\columnwidth]{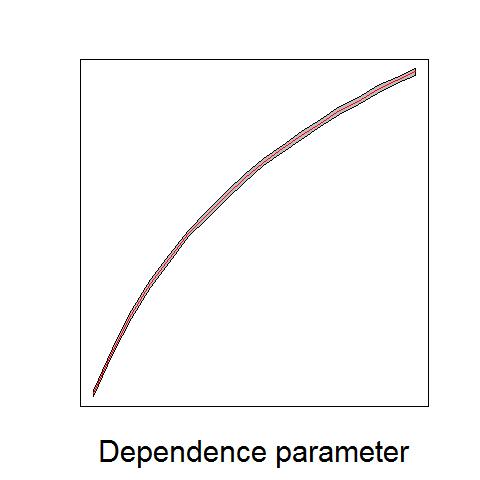}&\hspace{-0.17in} \includegraphics[trim=0.2in 0in 0.65in 0.65in,clip,height=0.45\columnwidth]{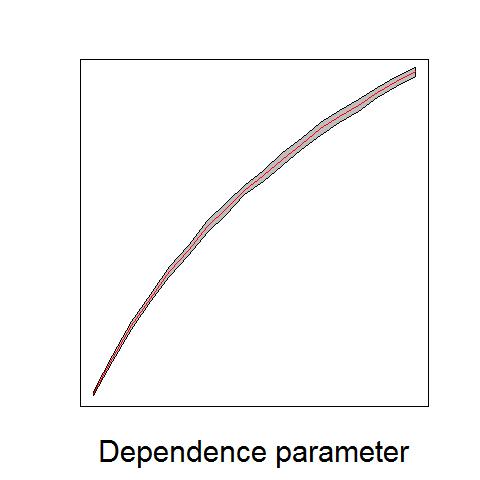}
\end{tabular}
\end{tabular}
\end{tabular}
\end{center}
\end{minipage}\hspace{0.04in}
\caption{Simulation study of the monotonicity of the empirical entropy in the copula dependence parameter. Shown
is the mean (red line) and 95\% range (black lines) of the empirical entropy over 50 random computations using
from $M=1000$ samples for each value of the dependence parameter (x-axis). We consider several popular bivariate (d=2) single parameter copula families (first and second rows), bivariate two parameters families (third row), as well as multivariate  (d=5) copula families (fourth row).}
\commentout{\mynote{Don't forget to improve figure by graying out range area instead of using 3 lines which are hard to distinguish}}
\label{BivariateSimulation}
\end{figure*}
\begin{figure*}
\begin{center}
\begin{tabular}{cc}
\begin{minipage}{0.45\columnwidth}
\caption{The impact of sample size on the monotonicity of the empirical entropy as a function of the dependence parameter of the generating copula. Shown is the mean monotonic trend measure (y-axis) as a function of the sample size (x-axis). The monotonicity is measured via the proportion of consecutive values of the dependence parameter, as defined in ~\eqref{EmpiricalMonotonicity}. Result shown is the average over 50 repetitions for the Clayton family. Results for other copula families are  qualitatively similar.}
\label{B}
\end{minipage} &
\begin{tabular}{c}
\includegraphics[trim=0.2in 0in 0.5in 0.5in,clip,height=0.5\columnwidth]{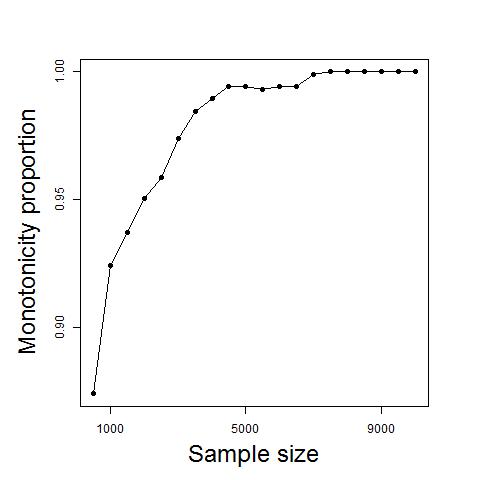}
\end{tabular}
\end{tabular}
\end{center}
\end{figure*}
 \commentout{Recently, various high-dimensional constructions that build on a collection of copulas have been suggested. Probably the most notable are the vine \cite{Bedford+Cooke,Kurowicka+Cooke:2002} and the Copula Bayesian Network (CBN) \cite{Kirshner,Elidan:2010}. Accordingly, one of the challenges that has been repeatedly tackled, both in the statistical and the machine learning communities is the problem of model selection. Importantly, as given a finite sample over N variables, the number of possible graph networks is exponential in $N$, examination of all possible networks is prohibited. Therefore most approaches rely on a greedy search procedure, where in each step, local modifications that maximize a score function, are performed.
 
 Somewhat interestingly, for both above mentioned constructions, the most common model selection methods (\cite{disSMann2013selecting,elidan2012lightning}, use Spearman's rho as the local score function and in each step, the algorithms stretch an edge between the two variables that share the largest magnitude of Spearman's rho. 
 
 While up to this point the justification for these approaches was the intuition, that higher value of Speraman's rho also lead to higher value of the log likelihood function, (which is the ultimate objective goal), our result, and in particular ~\lemref{MonotonicityOfSpearman}, confirms this intuition and supplies the asymptotical justification. 
 
 Note, however, that our results establish the monotonicity only asymptotically, while in practice, given a data set $\mathbb{D}$ that consists of N m-dimensional samples, and assume that the parametrization is given by a copula $C_{\theta}$, an essential assumption for these structure learning algorithms to succeed is that the empirical entropy is monotonic in the copula dependence parameter, that is:
 \begin{equation}
 \label{EmpiricalMonotonicity}
    \theta_1 \leq \theta_2 \Rightarrow \frac{\sum_{i=1}^N \log c_{\theta_1}(u_1[i],\ldots, u_m[i])}{N} \leq \frac{\sum_{i=1}^N \log c_{\theta_2}(u_1[i],\ldots, u_m[i])}{N}
 \end{equation}

 Therefore an interesting question is whether the monotonicity is evidential also in the finite sample case or alternatively how qualitative is the monotonicity in this case. 
 }

 The theory presented in this paper suggests that, in the case of an infinite number of samples generated from a given copula family, the entropy curve is monotonic in the copula dependence parameter. In practice, however, we almost always have access only to a finite number of samples. An obvious empirical question is thus whether monotonicity approximately holds given a reasonable number of samples, and what is the impact of the sample size on the entropy vs. dependence parameter curve.
 
 To answer these questions we explore the finite-sample behaviour of the copula entropy monotonicity via a simulation study. For each of the different theoretical scenarios discussed in the previous section, we choose representative popular copula families. We then generate $M=1000$ samples from the copula for different values of the dependence parameter. For each of these samples, we
 compute the entropy using the standard empirical
 estimator
 \[
     {\widehat H}(c_{\theta}) = -\frac{\sum_{{\bf u} \in \mathbb{D}}\log{c_{\theta}({\bf u})}}{M},
 \]
 where $\mathbb{D}$ denotes the set of $M$ samples generated from a copula $C_{\theta}$.

For each value of the dependence parameter, we repeat the above 50 times, and report the mean empirical entropy along with a $95\%$ confidence interval vs. the value of the dependence parameter. 
\figref{BivariateSimulation}, first two rows, show the results for single-parameter bivariate copula families. It is clear that with as little as $1000$ samples, near-monotonicity consistently holds for all families evaluated including both elliptical ones ({\bf Gaussian}, student-{\bf T}) and Archmedian copula families ({\bf Clayton}, {\bf Frank}, {\bf Gumbel}, {\bf Joe}, {\bf AMH}).

Next we turn to bivariate two parameters families. 
In this case we also need to test the monotonicity along the $\delta$ axis, since as formalized in  ~\corref{BivariatePrm}, monotonicity holds in $\theta$ for each \emph{fixed} value of $\delta$. Results are shown in \figref{BivariateSimulation}, third row,  for the {\bf BB1} and {\bf BB6} copula families, for two different values of $\delta$. Results for other values of $\delta$ as well as for the {\bf BB2} family were qualitatively similar. As before, near monotonicity is evident.     

To evaluate the finite-sample monotonicity in the multivariate case, we repeat the same evaluation for several copula families of dimension 5.
Results are shown in the last row of \figref{BivariateSimulation}, and are qualitatively similar for higher dimensions. As before, the empirical monotonicity is quiet impressive.

Finally, we evaluate the impact of the sample size on the extent to which the empirical entropy is monotonic in the copula dependence parameter. We repeat the same evaluation procedure for a wide range of sample sizes (from 500 to 10000).
To measure empirical monotonicity, for a given sample size, we measure the fraction of consecutive parameter values for which monotonicity holds. That is, we measure the empirical rate at which the following inequality holds:
\begin{equation}
\label{EmpiricalMonotonicity}
  \theta_1\leq \theta_2 \Rightarrow  \frac{\sum_{{\bf u} \in \mathbb{D}}\log{c_{\theta_1}({\bf u})}}{M}\leq \frac{\sum_{{\bf u} \in \mathbb{D}}\log{c_{\theta_2}({\bf u})}}{M},
\end{equation}
for two consecutive values of $\theta_1$ and $\theta_2$. We repeat this process 50 times and report the average. Results for the Clayton family
are shown in ~\figref{B}. As expected, with a greater sample size, monotonicity holds more frequently nearing 100\% at just 5000 training instances. Appealingly, even at much smaller sample sizes, monotonicity is quite appealing.

\commentout{\mynote{The curve has changed by the description in the text has not. Please verify that what we are saying is consistent with what you measured}}

%% file: Conclusions.tex
\section{Conclusion}
\label{conclusions}

In this work we establish a novel theoretical relationship between the main general purpose framework for capturing non-Gaussian real-valued distributions, namely copulas, and the ubiquitous Shannon's mutual information measure that is
used throughout the exact science to quantify the dependence between
random variables.

Our main result is that the mutual information between two variables (equivalently the copula entropy) is monotonic in the copula dependence parameter for a (very) wide range of copula families. Concretely, we
provide fairly general sufficient conditions for this monotonicity that cover
the vast majority of commonly used bivariate copulas, as well as a wide
class of multivariate Archimedean copulas.

As we demonstrated in our earlier work where the theory was substantially less developed \cite{tenzer2013speedy}, the monotonicity result also has practical merits. Specifically, it allows us to rank the entropy (equivalently the expected log-likelihood) of copula-based probabilistic models by
simple computation of association measures such as Spearman's rho.
This gives rise to highly efficient model selection in scenarios involving many interactions, e.g. when learning copula graphical models. The practical merit depends, of course, on the finite sample behavior of the entropy. Fortunately, as our simulations clearly show, near perfect monotonicity holds even with modest sample sizes.

Essentially all copula \emph{families} are constructed so as to span
some or all of the range between the independence and full dependence 
copula via the dependence parameter. We have not been able to identify
a single family where the above monotonicity does not hold empirically. However, monotonicity does appear to hold for the Plackett copula
familiy \cite{Nelsen} and identifying further sufficient conditions
remains a future challenge. Another direction if interesting is developing
finite sample theory that will explicitly quantify the amount by which
the empirical entropy can deviate from the expected monotonic behavior.